\documentclass[a4paper,leqno,10pt]{amsart}
\usepackage[top=3cm, bottom=3cm, left=3cm, right=3cm]{geometry}
\usepackage[T1]{fontenc}
\usepackage{lmodern}
\usepackage{amssymb}
\usepackage{amsmath}
\usepackage{mathrsfs}
\usepackage[usenames]{color}
\usepackage{bm}
\usepackage{setspace}
\usepackage{mathtools}
\allowdisplaybreaks

\begin{document}
%\renewcommand{\footnote}{}
%\pagenumbering{roman}
%\thispagestyle{empty}%\quad\newpage
%\thispagestyle{empty}
%\nonstopmode
%**************************************************************************

\newtheorem{theorem}{Theorem}[section] % Nummerierung
\newtheorem{proposition}[theorem]{Proposition}
\newtheorem{corollary}[theorem]{Corollary}
\newtheorem{lemma}[theorem]{Lemma}

\theoremstyle{definition}
\newtheorem{assumption}[theorem]{Assumption}
\newtheorem{definition}[theorem]{Definition}

\theoremstyle{definition} %%{remark}
\newtheorem{remark}[theorem]{Remark}
\newtheorem{remarks}[theorem]{Remarks}
\newtheorem{example}[theorem]{Example}
\newtheorem{examples}[theorem]{Examples}
%**************************************************************************
\newenvironment{pf}%
{\begin{sloppypar}\noindent{\bf Proof.}}%
{\hspace*{\fill}$\square$\vspace{6mm}\end{sloppypar}}
%**************************************************************************

\def\conj#1{{\overline{#1}}}
\def\A{{\mathbb{A}}}
\def\B{{\mathbb{B}}}
\def\R{{\mathbb{R}}}
\def\N{{\mathbb{N}}}
\def\Z{{\mathbb{Z}}}
\def\NO{{\mathbb{N}_0}}
\def\C{{\mathbb{C}}}
\def\M{{\mathbb{M}}}
\def\K{{\mathbb{K}}}
\def\E{{\mathbb{E}}}
\def\U{{\mathbb{U}}}
\def\F{{\mathbb{F}}}
\def\G{{\mathbb{G}}}
\def\I{{\mathbb{I}}}
\def\al{{\alpha}}
\def\bt{{\beta}}
\def\gm{{\gamma}}
\def\dl{{\delta}}
\def\ep{{\epsilon}}
\def\eps{{\varepsilon}}
\def\kp{{\kappa}}
\def\lm{{\lambda}}
\def\ph{{\phi}}
\def\rh{{\rho}}
\def\sg{{\sigma}}
\def\ek{{e_k}}
\def\ej{{e_j}}
\def\ei{{e_i}}
\def\ci{{c_i}}
\def\ck{{c_k}}
\def\cj{{c_j}}
\def\cit{{\tilde{c}_i}}
\def\ct{{\tilde{c}}}
\def\cnt#1{{\tilde{c}_#1}}
\def\cieq{{c_{i*}}}
\def\cjeq{{c_{j*}}}
\def\ckeq{{c_{k*}}}
\def\cneq#1{{c_{#1*}}}
\def\c{{c}}
\def\ct{{\tilde{c}}}
\def\ceq{{c_*}}
\def\ceqh{{\hat{c}_*}}
\def\cis{{c^\Blat_i}}
\def\cns#1{{c^\Blat_{#1}}}
\def\cn#1{{c_{#1}}}
\def\cks{{c^\Blat_k}}
\def\cjs{{c^\Blat_j}}
\def\cist{{\tilde{c}^\Blat_i}}
\def\cst{{\tilde{c}^\Blat}}
\def\cnst#1{{\tilde{c}^\Blat_#1}}
\def\ciseq{{c^\Blat_{i*}}}
\def\cjseq{{c^\Blat_{j*}}}
\def\ckseq{{c^\Blat_{k*}}}
\def\cnseq#1{{c^\Blat_{#1*}}}
\def\cs{c^\Blat}
\def\cst{{\tilde{c}^\Blat}}
\def\cseq{{c^\Blat_*}}
\def\cseqh{{\hat{c}^\Blat_*}}
\def\cio{{c_{i,0}}}
\def\ciot{{\tilde{c}_{i,0}}}
\def\cot{{\tilde{c}_{0}}}
\def\ciso{{c^\Sigma_{i,0}}}
\def\cisot{{\tilde{c}^\Sigma_{i,0}}}
\def\csot{{\tilde{c}^\Sigma_{0}}}
\def\co{{c_{0}}}
\def\cso{{c^\Sigma_{0}}}
\def\cs{{c^\Sigma}}
\def\di{{d_i}}
\def\dk#1{{d_{#1}}}
\def\dis{{d^\Sigma_i}}
\def\dks#1{{d^\Sigma_{#1}}}
\def\ris{{r^{\mathrm{sorp}}_i}}
\def\rs{{r^{\mathrm{sorp}}}}
\def\ric{{r^{\mathrm{ch}}_i}}
\def\rc{{r^{\mathrm{ch}}}}
\def\gin{{g^{\mathrm{in}}}}
\def\gini{{g^{\mathrm{in}}_i}}
\def\Om{{\Omega}}
\def\OmN{{\Om^N}}
\def\BOm{{\partial\Om}}
\def\Bin{{\Gamma_\mathrm{in}}}
\def\BinN{{\Gamma_\mathrm{in}^N}}
\def\Bnt{{\Gamma_\mathrm{nt}}}
\def\Bout{{\Gamma_\mathrm{out}}}
\def\BoutN{{\Gamma_\mathrm{out}^N}}
\def\Blat{{\Sigma}}
\def\BlatN{{\Blat^N}}
\def\BBlat{{\partial\Sigma}}
\def\BBlatN{{\partial\BBlat^N}}
\def\JT{{(0,T)}}
\def\Jinf{{(0,\infty)}}
\def\IT{{(0,\infty)}}
\def\Iinf{{[0, \infty)}}
\def\dt{{\partial_t}}
\def\dx{{\partial_x}}
\def\dxi{{\partial_{x_i}}}
\def\dn{{\partial_\nu}}
\def\dns{{\partial_{\nu_\Blat}}}
\def\lapl{{\Delta}}
\def\laplBelt{{\Delta_\Blat}}
\def\tlapl{{\tilde{\Delta}}}
\def\tlaplBelt{{\tlapl_\Blat}}
\def\kia{{k_i^\mathrm{ad}}}
\def\kna#1{{k_{#1}^\mathrm{ad}}}
\def\kid{{k_i^\mathrm{de}}}
\def\knd#1{{k_{#1}^\mathrm{de}}}
\def\Kad{{K^\mathrm{ad}}}
\def\Kde{{K^\mathrm{de}}}
\def\Kab#1#2{{K^{#1}_{#2}\,}}
\def\sqrtKad{{\left(\Kad\right)^{1/2}}}
\def\sqrtKde{{\left(\Kde\right)^{1/2}}}
\def\pwrKad#1{{\left(\Kad\right)^{#1}}}
\def\pwrKde#1{{\left(\Kde\right)^{#1}}}
\def\matin{{\mathrm{in}}}
\def\maton{{\mathrm{on}}}
\def\LebSpc#1{{L_{#1}}}
\def\SobSpc#1#2{{W^{#1,#2}}}
\def\HilbSpc#1{{H^{#1}}}
\def\SobSlobSpc#1#2{{W^{#1}_{#2}}}
\def\LebSpcloc#1{{L_{#1, loc}}}
\def\SobSpcloc#1#2{{W^{#1,#2}_{loc}}}
\def\HilbSpcloc#1{{H^{#1}_{loc}}}
\def\Lp{{\LebSpc{p}}}
\def\Wkp{{\SobSpc{k}{p}}}
\def\LOm{{\LebSpc{2}(\Om)}}
\def\LOmN{{\LebSpc{2}(\Om)^N}}
\def\LpOmN#1{{\LebSpc{#1}(\Om)^N}}
\def\LOmNN{{\LebSpc{2}(\Om)^{N \times N}}}
\def\LBlat{{\LebSpc{2}(\Blat)}}
\def\LBlatN{{\LebSpc{2}(\Blat)^N}}
\def\LpBlatN#1{{\LebSpc{#1}(\Blat)^N}}
\def\LBlatNN{{\LebSpc{2}(\Blat)^{N \times N}}}
\def\LOmBlat{{\LOm \times \LBlat}}
\def\LOmBlatN{{\LOmN \times \LBlatN}}
\def\LpOmBlatN#1{{\LpOmN{#1} \times \LpBlatN{#1}}}
\def\LOmBlatNN{{\LOmNN \times \LBlatNN}}
\def\UO#1{{\U_{#1}^\Om}}
\def\EO#1{{\E_{#1}^\Om}}
\def\ES#1{{\E_{#1}^\Blat}}
\def\FO#1{{\F_{#1}^\Om}}
\def\FS#1{{\F_{#1}^\Blat}}
\def\FOS#1{{\F_{#1}^{\Om, \Blat}}}
\def\FOSI#1{{\F_{{#1}, I}^{\Om, \Blat}}}
\def\Gi#1{{\G_{#1}^\mathrm{in}}}
\def\Go#1{{\G_{#1}^\mathrm{out}}}
\def\GS#1{{\G_{#1}^\Blat}}
\def\UOp{{\UO{p}}}
\def\EOp{{\EO{p}}}
\def\ESp{{\ES{p}}}
\def\FOp{{\FO{p}}}
\def\FSp{{\FS{p}}}
\def\FOSp{{\FOS{p}}}
\def\FOSpI{{\FOSI{p}}}
\def\Gip{{\Gi{p}}}
\def\Gop{{\Go{p}}}
\def\GSp{{\GS{p}}}
\def\InitSpc#1#2{{\I_{#1}(#2)}}
\def\IOp{{\InitSpc{p}{\Om}}}
\def\ISp{{\InitSpc{p}{\Blat}}}
\def\trcrega{{1/2 - 1/{2p}}}
\def\trcregb{{1 - 1/p}}
\def\trcregc{{2 - 2/p}}
\def\Cp{{C_p}}
\def\CT{{C_T}}
\def\Cv{{C_v}}
\def\CuO{{C_{u, \Om}}}
\def\Cpq{{C_p^2}}
\def\CTq{{C_T^2}}
\def\Cvq{{C_v^2}}
\def\CuOq{{C_{u, \Om}^2}}
\def\Av{{$\textup{A}^\textup{vel}$}}
\def\Avi{{$\textup{A}^\textup{vel}_\textup{in}$}}
\def\Avo{{$\textup{A}^\textup{vel}_\textup{out}$}}
\def\AFs{{$\textup{A}_\textup{F}^\textup{sorp}$}}
\def\AMs{{$\textup{A}_\textup{M}^\textup{sorp}$}}
\def\ABs{{$\textup{A}_\textup{B}^\textup{sorp}$}}
\def\AFc{{$\textup{A}_\textup{F}^\textup{ch}$}}
\def\ANc{{$\textup{A}_\textup{N}^\textup{ch}$}}
\def\APc{{$\textup{A}_\textup{P}^\textup{ch}$}}
\def\APe{{$\textup{A}_\textup{P}^\textup{eq}$}}
\def\ARe{{$\textup{A}_\textup{R}^\textup{eq}$}}
\def\AIe{{$\textup{A}_\textup{I}^\textup{eq}$}}

% Equation Spacing
\newcommand{\shrinkdisplayskips}  {\addtolength{\abovedisplayskip}{-0.5em}\addtolength{\abovedisplayshortskip}{-0.5em}\addtolength{\belowdisplayskip}{-0.5em}\addtolength{\belowdisplayshortskip}{-0.5em}}
\newcommand{\unshrinkdisplayskips}{\addtolength{\abovedisplayskip}{ 0.5em}\addtolength{\abovedisplayshortskip}{ 0.5em}\addtolength{\belowdisplayskip}{ 0.5em}\addtolength{\belowdisplayshortskip}{ 0.5em}}
\unshrinkdisplayskips

% Array Arrangements
\newcommand{\narrowarray}{\setlength{\arraycolsep}{0.2em}}
\newcommand{\normalarray}{\setlength{\arraycolsep}{0.4em}}
\newcommand{\widearray}  {\setlength{\arraycolsep}{0.6em}}
\narrowarray

\DeclarePairedDelimiter{\tmpnorm}{\lVert}{\rVert}
\DeclarePairedDelimiter{\tmpabs}{\lvert}{\rvert}
\DeclarePairedDelimiter{\tmpsclprd}{(}{)}
\makeatletter
\newcommand{\@normstar}[2]{\tmpnorm*{#2}_{#1}}
\newcommand{\@normnostar}[3][]{\tmpnorm[#1]{#3}_{#2}}
\newcommand{\norm}{\@ifstar\@normstar\@normnostar}
\newcommand{\@normQuadstar}[2]{\tmpnorm*{#2}_{#1}^2}
\newcommand{\@normQuadnostar}[3][]{\tmpnorm[#1]{#3}_{#2}^2}
\newcommand{\normQuad}{\@ifstar\@normQuadstar\@normQuadnostar}
\newcommand{\@absstar}[1]{\tmpabs*{#1}}
\newcommand{\@absnostar}[2][]{\tmpabs[#1]{#2}}
\newcommand{\abs}{\@ifstar\@absstar\@absnostar}
\newcommand{\@sclprdstar}[3]{\tmpsclprd*{#2, #3}_{#1}}
\newcommand{\@sclprdnostar}[4][]{\tmpsclprd[#1]{#3, #4}_{#2}}
\newcommand{\sclprd}{\@ifstar\@sclprdstar\@sclprdnostar}
\makeatother

\def\Ball#1#2#3{{\B_{#1}^{#2}\left( #3 \right)}}
\def\BallCompl#1#2#3{{\overline{\B}_{#1}^{#2}\left( #3 \right)}}
\def\hc{{Heterogeneous Catalysis}\space}
\def\HC{{HETEROGENEOUS CATALYSIS}\space}
\def\id{{\textrm{Id}}}
\def\diag{{\textrm{diag}}}

\def\Poinc{{Poincar\'{e}}}
\def\Frech{{Fr\'{e}chet}}
\def\PoincS{{\Poinc\space}}
\def\FrechS{{\Frech\space}}

\hyphenation{Lipschitz}

%**************************************************************************
\sloppy
%**************************************************************************
\title[Stability Analysis for a class of \HC models]
	{Stability Analysis for a class of \\ \HC models}

\author[C.\ Gesse]{Christian Gesse}
\address{\noindent
	Mathematisches Institut -- Heinrich-Heine-Uni\-ver\-sit\"at D\"usseldorf \newline\indent
	Universit\"atsstr.\ 1, 40225 D\"usseldorf, Germany}
\email{christian.gesse@hhu.de}

\author[M.\ K\"ohne]{Matthias K\"ohne}
\address{
	Mathematisches Institut -- Heinrich-Heine-Uni\-ver\-sit\"at D\"usseldorf \newline\indent
	Universit\"atsstr.\ 1, 40225 D\"usseldorf, Germany}
\email{matthias.koehne@hhu.de}

\author[J.\ Saal]{J\"urgen Saal}
\address{
	Mathematisches Institut -- Heinrich-Heine-Uni\-ver\-sit\"at D\"usseldorf \newline\indent
	Universit\"atsstr.\ 1, 40225 D\"usseldorf, Germany}
\email{juergen.saal@hhu.de}



\date{\today}

\parskip0.5ex plus 0.5ex minus 0.5ex

\begin{abstract}
We prove stability for a class of heterogeneous catalysis models 
in the $L_p$-setting.
We consider a setting in a finite three-dimensional pore of 
cylinder-like geometry, with the lateral walls acting as a catalytic surface.
Under a reasonable condition on the involved parameters, we show that given 
equilibria are normally stable, i.e.\ solutions are
attracted at an exponential rate. The potential incidence of instability 
is discussed as well.
\end{abstract}
\maketitle
%%%%%%%%%%%%%%%%%%%%%%%%%%%%%%%%%%%%

%{\bf 2020 Mathematics Subject Classification.}
%Primary: 35K57; Secondary: 35K55, 35R01, 80A32.
 
%\tableofcontents

%\doublespacing

%%%%%%%%%%%%%%%%%%%%%%%%%%%%%%%%%%%%%%%%%%%%%%%%%%%%%%%%%%%%%%%%%%%%%%%%%%%%%
\section{Introduction}

As a key technology in chemical engineering, besides for the increase of the speed of chemical reactions, catalysis is also employed to change the selectivity in favor of a desired product against other possible output components of a chemical reaction network.
One can differentiate between {\itshape homogeneous catalysis}, where the catalyst itself is in the same phase as the other reactants, and {\itshape heterogeneous catalysis}, where the catalyst is present in a different phase, usually a solid wall.
The latter case is advantageous concerning the separation of the products from the catalytic material.
However, a high area-to-volume ratio is required, which is given, for instance, in case of the walls of a porous medium.
For more information on heterogeneous catalysis we refer to \cite{levenspiel:chemical_reaction_eng,aris:mathematical_theory,white:het_cat}.
In this regard, this note is devoted to the study of the stability of the 
following prototype model for heterogeneous catalysis in a cylindric domain:

Let $A \subset \R^2$ be a bounded, simply connected $C^2$-domain and let $h > 0$.
Let $\Om := A \times (0, h)$ be a finite
three-dimensional cylinder. We decompose the smooth part of the boundary
$\partial\Om$ into the inflow surface $\Bin = A \times \{ 0 \}$, the
outflow surface $\Bout = A \times \{ h \}$ and the lateral surface
$\Blat = \partial A \times (0,h)$. We consider the following system of balance equations:
\begin{equation}\label{het_cat_sys}
	\begin{array}{rclll}
		\dt\ci + (u \cdot \nabla)\ci - \di\lapl\ci &=& 0 &\qquad\matin &\JT\times\Om,\\[0.5em]
		\dt\cis - \dis\laplBelt\cis &=& \ris(\ci, \cis) + \ric(\cs) &\qquad\maton &\JT\times\Blat,\\[0.5em]
		(u \cdot\nu)\ci - \di\dn\ci &=& \gini &\qquad\maton &\JT\times\Bin,\\[0.5em]
		-\di\dn\ci &=& \ris(\ci, \cis) &\qquad\maton \, &\JT\times\Blat,\\[0.5em]
		-\di\dn\ci &=& 0 &\qquad\maton &\JT\times\Bout,\\[0.5em]
		-\dis\dns\cis &=& 0 &\qquad\maton &\JT\times\BBlat,\\[0.5em]
		\ci\vert_{t=0} &=& \cio &\qquad\matin &\Om,\\[0.5em]
		\cis\vert_{t=0} &=& \ciso &\qquad\maton &\Blat,
		\end{array}
\end{equation}
where $\c := (\ci)_{i=1}^N$ denote the bulk concentrations and $\cs := (\cis)_{i=1}^N$ denote the surface concentrations of the involved species $(C_i)_{i=1}^N$.
The velocity field $u$ is assumed to be given (and sufficiently smooth).
\pagebreak

In this paper, we show stability of positive equilibria $(\ceq, \cseq)$ for (\ref{het_cat_sys}) in the $\Lp$-setting.
We restrict our choice for the {\itshape sorption rates} $\ris$ to the linear case
\begin{equation}\label{sorp_func}
	\ris(\ci, \cis) = \kia\ci - \kid\cis,
\end{equation}
where $\kia, \kid > 0$.
For the {\itshape chemical reaction rates} $\ric$ we assume that the reaction of $N$ species is given as a reversible reaction
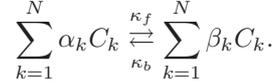
\begin{displaymath}
	\sum_{k = 1}^{N} \al_k C_k \overset{\kp_f}{\underset{\kp_b}{\rightleftarrows}} \sum_{k = 1}^{N} \bt_k C_k.
\end{displaymath}
Here, $\kp_f > 0$ denotes the forward reaction rate and $\kp_b > 0$ the backward reaction rate, while $(\al_k)_{k=1}^N,\ (\bt_k)_{k=1}^N \in (\{ 0 \} \cup [1, \infty))^N \setminus \{ 0 \}^N$ denote the stoichiometric coefficients.
The reaction rate for this reaction is then given as
\begin{equation}\label{chem_func}
	\ric(\cs) := (\al_i - \bt_i)\left(\kp_b \prod_{k=1}^{N} (\cks)^{\bt_k} - \kp_f \prod_{k=1}^{N} (\cks)^{\al_k}\right).
\end{equation}

In the model described above the educt species are transported from the bulk phase to the 
lateral surface $\Sigma$, where they adsorb with rate $\kia\ci$,
see \eqref{sorp_func}. The adsorbed molecules 
react with rate $\ric$ with other adsorbed molecules.
The product molecules are desorbed to the bulk again with the rate $\kid\cis$.
Note that all chemical reactions take place on the lateral surface $\Sigma$,
where the catalyst is assumed to be present.

The equations modeling heterogeneous catalysis processes considered in this work have been proposed in \cite{maier:het_cat}, where a mathematical analysis of linear and nonlinear local well-posedness is carried out.
Additionally, global well-posedness is proved under the additional assumption of a triangular structure of the chemical reaction rate.
For further details on the model we refer
to \cite{maier:het_cat} and the references therein. There are more
recent results on the mathematical modeling of the heterogeneous
catalysis process, see e.g. \cite{soucek:het_cat} for a detailed
approach modeling a coupled system of equations in a suitable
thermodynamic framework. In \cite{augner:fast_sorption} several limit
models are derived taking into account the time scale on which the
chemical reaction and the sorption occur. Additionally, a three
component model problem is analyzed in terms of well-posedness,
positivity of solutions, blow-up criteria and a-priori bounds. The
approach is extended to more general systems in \cite{augner:het_cat}.
Recent results regarding global well-posedness of volume-surface
reaction-diffusion systems are presented in
\cite{morgan:global_well_posedness}. However, no cylindrical structure
with inflow- and outflow surface is considered in these works and the
results do not cover stability or instability of equilibria. In
\cite{schnaubelt:dynamic_nonlinear_bc} a general theory regarding
stable and unstable manifolds is developed for quasilinear problems with
nonlinear dynamical boundary conditions. However, the functional
analytic setting there differs from
the one we use in this work as in \cite{schnaubelt:dynamic_nonlinear_bc} 
the quantities on the boundary are functions in
the trace space of the space for the quantities in the bulk. For system
\eqref{het_cat_sys} we find it more appropriate to have the same 
regularity for both, the quantities in the bulk and on the boundary. 
Thus, as it seems, the two appoaches are not 
comparable.

This note is organized as follows: In Section~\ref{sec_not} the notation
is settled. In Section~\ref{sec_mr} we recall a result on maximal
regularity for a linearization of \eqref{het_cat_sys} derived in
\cite{maier:het_cat}. Based on this and the 
principle of linearized stability, see \cite{pruess:conv,pruess:movint},
we prove the main result of this paper on stability of equilibria for system
\eqref{het_cat_sys} in Section~\ref{sec_stab}. Conditions that yield instabilities
are discussed in Section~\ref{sec_instab}.  
The paper ends with a concluding debate of the obtained results
in Section~\ref{sec_concl}.

%%%%%%%%%%%%%%%%%%%%%%%%%%%%%%%%%%%%%%%%%%%%%%%%%%%%%%%%%%%%%%%%%%%%%%%%%%%%%
\section{Notation}\label{sec_not}
Throughout this paper, for $\K \in \{\R, \C \}$ and $n \in \N$ let $|\,\!\cdot\,\!|$ denote the euclidean norm on $\K^n$ and let $|\,\!\cdot\,\!|_2$ denote the induced spectral norm on $\K^{n \times n}$.
Furthermore, let $\SobSpc{k}{p}(G,X)$ denote the the $X$-valued Sobolev space and $\SobSlobSpc{s}{p}(G,X)$ denote the $X$-valued Sobolev-Slobodeckij space for a Banach space $X$ (with norm $\norm{X}{\cdot}$),
a domain $G \subset \R^n$, $k \in \NO := \N \cup \{0\} $ and $s \in (0, \infty) \setminus \N$. Consequently, $\Lp(G) := \SobSpc{0}{p}(G)$ denotes the corresponding Lebesgue spaces.
Sometimes we will omit $X$ and write $\SobSpc{k}{p}(G)$, $\SobSlobSpc{s}{p}(G)$ and $\Lp(G)$ if no confusion is likely to arise.
Moreover, $\HilbSpc{k}(G) := \SobSpc{k}{2}(G)$. Additionally, $\sclprd{H}{\cdot}{\cdot}$ denotes the scalar product on a Hilbert space $H$.

In the case of $\LebSpc{2}(G)^n$ we have
\begin{align*}
	\sclprd{\LebSpc{2}(G)^n}{u}{v} = \int_G \sclprd{\K^n}{u}{v} \, dx, \qquad u,v \in \LebSpc{2}(\G)^n
\end{align*}
as well as
\begin{align*}
	\norm{\LebSpc{2}(G)^n}{u}^2 = \int_G |u|^2 \, dx, \qquad u \in \LebSpc{2}(\G)^n.
\end{align*}
If $Y((0, T), X)$ constitutes a Banach function space of $X$-valued functions on $(0, T)$,
we write ${}_0 Y((0, T), X)$ for the subspace of functions with trace zero at $t = 0$, provided that this trace is well-defined.

For a domain $G \subset \R^n$ we denote by $\partial G$ its boundary and by $\nu := \nu(x)$ the outer unit normal vector at $x \in \partial G$.
An open ball in $X$ with radius $r > 0$ and center $x \in X$ is denoted as $\Ball{r}{X}{x}$.
For $\al = (\al_1,\dots,\al_n) \in [0, \infty)^n$ and $x = (x_1,\dots,x_n) \in \K^n$ we set
\begin{align*}
	x^\al := \prod_{k=1}^{n} x_k^{\al_k}.
\end{align*}
Moreover, we write
\begin{align*}
	\diag(d_1, ... , d_n) := \begin{pmatrix}
	d_1 & &\\
	& \ddots & \\
	& & d_n
	\end{pmatrix}
\end{align*}
for a diagonal matrix in $\K^{n \times n}$ with diagonal values $d_1, \dots , d_n \in \K$.

If $\lm \in \C$, then $\Re \lm$ denotes the real part and $\Im \lm$ denotes the imaginary part of $\lambda$, respectively.
Additionally, $\C_\pm := \{ \lm \in \C \, : \, \pm\Re\lm > 0\}$.

For a closed operator $A: D(A) \subset X \rightarrow X$ let $\sg(A)$ denote the spectrum and $\rh(A)$ denote the resolvent set of $A$.

Let $\partial_j f := \partial_{x_j} f$ be the partial derivative with respect to the variable $x_j$ of a (sufficiently smooth) mapping $f: U \subset \K^n \rightarrow X$.
Then we have
\begin{align*}
	\nabla f := (\partial_1 f, \dots , \partial_n f)^T
\end{align*}
as the gradient of $f$ and
\begin{align*}
	\lapl f := \sum_{j=1}^n \partial_j^2 f
\end{align*}
as the Laplacian of $f$.
If a mapping $\phi:M \rightarrow X$ is defined on the (sufficiently smooth) surface $M$ and can be extended to a (one-sided) neighborhood of $M$,
then $\nabla_M  \phi:=  (\nabla \phi)\vert_M - \nu (\nu \cdot (\nabla \phi)\vert_M)$ denotes the surface gradient on $M$ whereas $\lapl_M \phi := \nabla_M \cdot \nabla_M \phi$ denotes the Laplace-Beltrami operator on $M$.
%%%%%%%%%%%%%%%%%%%%%%%%%%%%%%%%%%%%%%%%%%%%%%%%%%%%%%%%%%%%%%%%%%%%%%%%%%%%%
\section{Maximal regularity}\label{sec_mr}
The model (\ref{het_cat_sys}) considered in this paper was originally considered by Bothe, K\"ohne, Maier and Saal in \cite{maier:het_cat}.
Besides local and global existence with some restrictions regarding the sorption rates $\ris$ and the chemical reaction rates $\ric$ they showed maximal regularity for a linearization of (\ref{het_cat_sys}).
We will state their result of maximal regularity here, as it plays a crucial role for our stability analysis.
First, we define the corresponding solution spaces for $\ci$ and $\cis$, which are
\begin{align*}
	\EOp(T) &:= \SobSpc{1}{p}(\JT, \Lp(\Om)) \cap \Lp(\JT, \SobSpc{2}{p}(\Omega)),\\
	\ESp(T) &:= \SobSpc{1}{p}(\JT, \Lp(\Blat)) \cap \Lp(\JT, \SobSpc{2}{p}(\Blat)).
\end{align*}
The data spaces are derived using appropriate trace theorems, which leads to the spaces
\begin{align*}
	\FOp(T) &:= \Lp(\JT \times \Om),\\
	\FSp(T) &:= \Lp(\JT \times \Blat),\\
	\Gip(T) &:= \SobSlobSpc{\trcrega}{p}(\JT, \Lp(\Bin)) \cap \Lp(\JT, \SobSlobSpc{\trcregb}{p}(\Bin)),\\
	\GSp(T) &:= \SobSlobSpc{\trcrega}{p}(\JT, \Lp(\Blat)) \cap \Lp(\JT, \SobSlobSpc{\trcregb}{p}(\Blat)),\\
	\Gop(T) &:= \SobSlobSpc{\trcrega}{p}(\JT, \Lp(\Bout)) \cap \Lp(\JT, \SobSlobSpc{\trcregb}{p}(\Bout)),\\
	\IOp &:= \SobSlobSpc{\trcregc}{p}(\Om),\\
	\ISp &:= \SobSlobSpc{\trcregc}{p}(\Blat).
\end{align*}
Now, we can define the data space for the right-hand side of (\ref{het_cat_sys}) as
\begin{displaymath}
	\FOSp(T) := \FOp(T) \times \FSp(T) \times \Gip(T) \times \GSp(T) \times \Gop \times \{0\}
\end{displaymath}
and the corresponding space including the initial data as
\begin{displaymath}
	\FOSpI(T) := \FOSp(T) \times \IOp \times \ISp.
\end{displaymath}
Additionally, we impose the following restrictions regarding the velocity field $u$.
\begin{itemize}
	\item\textbf{(\Av)} Let $u$ denote a given velocity-field of regularity 
	\begin{displaymath}
	u \in \UOp(T) := \SobSpc{1}{p}(\JT, \Lp(\Om, \R^3)) \cap \Lp(\JT, \SobSpc{2}{p}(\Om, \R^3))
	\end{displaymath}
	fulfilling
	\begin{displaymath}\label{vel_u_assumptions}
	u \cdot \nu \leqslant 0 \ \maton \ \Bin, \quad u \cdot \nu = 0 \ \maton \ \Blat, \quad u \cdot \nu \geqslant 0 \ \maton \ \Bout
	\end{displaymath}
	and $\nabla \cdot u = 0$ in $\Omega$.
\end{itemize}
Now, the principal linearization of (\ref{het_cat_sys}) is given as
\begin{equation}\label{het_cat_sys_maxreg}
	\begin{array}{rclll}
		\dt\ci + (u \cdot \nabla)\ci - \di\lapl\ci &=& f_i & \quad \matin \, &\JT\times\Om,\\[0.5em]
		\dt\cis - \dis\laplBelt\cis &=& f_i^\Blat & \quad \maton \, &\JT\times\Blat,\\[0.5em]
		(u \cdot\nu)\ci - \di\dn\ci &=& \gini & \quad \maton \, &\JT\times\Bin,\\[0.5em]
		-\di\dn\ci &=& g_i^\Blat & \quad \maton \, &\JT\times\Blat,\\[0.5em]
		-\di\dn\ci &=& g_i^{\textrm{out}} & \quad \maton \, &\JT\times\Bout,\\[0.5em]
		-\dis\dns\cis &=& 0 & \quad \maton \, &\JT\times\BBlat,\\[0.5em]
		\ci\vert_{t=0} &=& \cio & \quad \matin \, &\Om,\\[0.5em]
		\cis\vert_{t=0} &=& \ciso & \quad \maton \, &\Blat.
	\end{array}
\end{equation}
The result for maximal regularity in the $\LebSpc{p}$-setting is given in \cite{maier:het_cat} and reads as follows.
\begin{theorem}\label{thm:maxreg}{\upshape (\cite[Prop.\ 4.1]{maier:het_cat}).}
	Let $T > 0$, let $J = (0, T) \subset \R$ and let $\frac{5}{3} < p < \infty$ with $p \neq 3$.
	Suppose the velocity field $u$ satisfies (\Av). Then (\ref{het_cat_sys_maxreg}) admits a unique solution
	\begin{align*}
		(\ci, \cis) \in \EOp(T) \times \ESp(T),
	\end{align*}
	if and only if the data satisfy the regularity condition
	\begin{align*}
		(f_i, f_i^\Blat, \gini, g_i^\Blat, g_i^{\textrm{out}}, 0, \cio, \ciso) \in \FOSpI(T)
	\end{align*}
	and, in the case $p > 3$, the compatibility conditions
	\begin{align}\label{compatibility_condition}
		\begin{aligned}
			\cio u(0) \cdot \nu - \di\dn\cio &= \gini(0) &\maton\; &\Bin,\\
			-\di\dn\cio &= \ris(\cio, \ciso) &\maton\; &\Blat,\\
			-\di\dn\cio &= 0 &\maton\; &\Bout,\\
			-\dis\dns\ciso &= 0 &\maton\; &\BBlat.
		\end{aligned}
	\end{align}
	Additionally, the corresponding solution operator ${}_0\mathcal{S}_T$ w.r.t.\ homogeneous initial conditions satisfies
	\begin{align*}
		\norm{\mathscr{L}({}_0 \FOSp(\tau)^N, {}_0\EOp(\tau)^N \times {}_0\ESp(\tau)^N)}{{}_0 \mathcal{S}_T} \leqslant M, \qquad 0 < \tau < T,
	\end{align*}
	for a constant $M > 0$ that is independent of $0 < \tau < T$.
\end{theorem}
%%%%%%%%%%%%%%%%%%%%%%%%%%%%%%%%%%%%%%%%%%%%%%%%%%%%%%%%%%%%%%%%%%%%%%%%%%%%
\section{Stability in the $\LebSpc{p}$-setting}\label{sec_stab}
In this section we prove stability for (\ref{het_cat_sys}) in the $\LebSpc{p}$-setting for $p \in [2, \infty) \setminus \{3\}$ and for a fixed sorption rate, but with a wide choice of reactions and equilibria. 
\begin{remark}\label{eq_chem_bal}
	As an example of an equilibrium one may choose the equilibrium $(\cieq, \ciseq)_{i=1,\dots,N}$ of (\ref{het_cat_sys}) as the constant equilibrium of chemical balance, i.e.
	\begin{align}\label{chem_balance_1}
	\cieq \equiv \psi_i > 0, \quad \ciseq \equiv \xi_i > 0, \qquad i = 1,\dots,N,
	\end{align}
	where
	\begin{align}\label{chem_balance_2}
	\psi_i = \frac{\kid}{\kia} \xi_i \quad \textrm{and} \quad \kp_b \prod_{k=1}^{N} (\xi_k)^{\bt_k} = \kp_f \prod_{k=1}^{N} (\xi_k)^{\al_k} ,
	\end{align}
	which ensures that $r_i(\cseq) = 0$.
	Here, we have to assume that the inflow profile fulfills $\gini \leqslant 0$ and $\gini \neq 0$ on $\Bin$ for $i =1,\dots,N$.
	Now, if the velocity profile at the inflow satisfies
	\begin{align*}
	(u \cdot \nu) = \frac{\kia}{\kid \xi_i} \gini \quad \maton \ \Blat, \qquad i = 1, \dots, N,
	\end{align*}
	then (\ref{chem_balance_1}) and (\ref{chem_balance_2}) ensure that $(\cieq, \ciseq)_{i=1,\dots,N}$ constitutes and equilibrium of (\ref{het_cat_sys}).
\end{remark}
This example motivates the following conditions, which we will assume to be fulfilled for any equilibrium $(\ceq, \cseq)$ we may choose in order to show stability:
\begin{itemize}
	\item\textbf{(\APe)} The equilibrium is non-negative, i.e.
	\begin{displaymath}
	\cieq \geqslant 0 \ \matin \ \Om, \quad \ciseq \geqslant 0 \ \maton \ \Blat, \qquad i = 1,\dots,N.
	\end{displaymath}
	\item\textbf{(\ARe)} The equilibrium fulfills the following regularity conditions:
	\begin{displaymath}
	\cieq \in \SobSpc{2}{p}(\Om), \quad \ciseq \in \SobSpc{2}{p}(\Blat).
	\end{displaymath}
	\item\textbf{(\AIe)} The equilibrium is isolated,
	i.e.\ $\Ball{\varepsilon}{\SobSpc{2}{p}(\Om) \times \SobSpc{2}{p}(\Blat)}{(\ceq, \cseq)} \setminus \{(\ceq, \cseq)\}$ does not contain another equilibrium for some $\varepsilon > 0$.
\end{itemize}
Furthermore, we impose an additional condition regarding the velocity field $u$:
\begin{itemize}
	\item\textbf{(\Avi)} The velocity field has non-trivial inflow, i.e.\ $u \cdot \nu \neq 0$ on $\Bin$.
\end{itemize}

Before we state our main result, we recall the following version of Poincar{\'e}'s inequality.
\begin{remark}\label{poincare}
	By \cite[Lemma 10.2 (vi)]{tartar:SobSpc} we have the following: Let $\varnothing \neq M \subset \R^n$ be a open and let $1 \leqslant p \leqslant \infty$.
	Let $V \subset \SobSpc{1}{p}(M)$ be a linear subspace. If the injection $V \hookrightarrow \Lp(M)$ is compact and the constant function $u \equiv 1$ does not belong to $V$, then there exists a constant  $C > 0$ s.t.
	\begin{align*}
		\norm{\Lp(M)}{u} \leqslant C \norm{\Lp(M)}{\nabla u}, \qquad u \in V,
	\end{align*}
	and one says that the \PoincS inequality holds. This assertion also holds if $M$ is replaced by the lateral boundary $\Blat$ of a cylindrical domain $\Om = A \times (0,h)$ with a simply connected $C^2$-domain $A \subset \R^2$.
\end{remark}

Now, our main result reads as follows.
\begin{theorem}\label{thm: stability}{\upshape (Stability in $\LebSpc{p}$).}
	Let $p \in [2, \infty) \setminus \{3\}$, let $T = \infty$ and let $\gini \in \Gip$ for $i = 1,\dots,N$. Let the sorption rates be given as 
	\begin{displaymath}
	\ris(\ci, \cis) := \kia \ci - \kid \cis, \qquad i = 1,...,N,
	\end{displaymath}
	and the reaction rates as
	\begin{displaymath}
	\ric(\cs) := (\al_i - \bt_i)\left(\kp_b \left(\cs\right)^\bt - \kp_f \left(\cs\right)^\al\right), \qquad i = 1,...,N,
	\end{displaymath}
	with $\kia, \kid > 0$, $\kp_b, \kp_f > 0$ and $\al, \bt \in (\{ 0 \} \cup [1, \infty))^N \setminus \{ 0 \}^N$.
	Assume that $(\ceq, \cseq)$ is an equilibrium of (\ref{het_cat_sys}) satisfying (\APe), (\ARe), (\AIe) and that the velocity field $u$ satisfies the additional condition (\Avi).
	Let
	\begin{align*}
	\max\nolimits_{\Blat} \abs{a}\abs{b(\cseq)} \leqslant \frac{1}{C_P},
	\end{align*}
	where $C_P > 0$ denotes the \PoincS constant on $\Blat$ (cf. Remark \ref{poincare}) and 
	\begin{align*}
		\left. \begin{aligned}
		a_k &:= (\al_k - \bt_k),\;\\ 
		b_k &:= b_k (\cseq) := \left(\kp_b \bt_k \left( \cseq \right)^{\bt - \ek} - \kp_f \al_k \left( \cseq \right)^{\al - \ek}\right)
		\end{aligned}\right.	
	\end{align*}
	for $k=1,\dots,N$.
	Then there exists $\rh > 0$ s.t.\ for 
	\begin{align*}
		(\co, \cso) \in \Ball{\rho}{\IOp^N \times \ISp^N}{(\ceq, \cseq)},
	\end{align*}
	where in case $p > 3$ the conditions (\ref{compatibility_condition}) have to be fulfilled, there exists a unique global solution $(\c, \cs)$ to (\ref{het_cat_sys}) satisfying
	\begin{align*}
	(\c, \cs) \in \SobSpcloc{1}{p}(\R_+, \LebSpc{p}(\Om)^N \times \LebSpc{p}(\Blat)^N) \cap \LebSpcloc{p}(\R_+, \SobSpc{2}{p}(\Om)^N \times \SobSpc{2}{p}(\Blat)^N).
	\end{align*}
	Moreover, the equilibrium $(\ceq, \cseq)$ is exponentially stable in $\IOp^N \times \ISp^N$.
\end{theorem}

\textit{Proof:} We want to apply the principle of linearized stability (cf. \cite{pruess:conv,pruess:movint}).
In order to shorten the notation, we write e.g.\ $c = (c_1,...,c_N)^T$ and, similarly, for all other appearing quantities.

Let $(\ceq, \cseq) \in \SobSpc{2}{p}(\Om)^N \times \SobSpc{2}{p}(\Blat)^N$ be an equilibrium fulfilling the assumptions (\APe), (\ARe) and (\AIe).
Moreover, assume that the velocity field $u$ satisfies the additional condition (\Avi).
We will proceed in three steps.

\textbf{Step 1: Translation of the system and mapping properties.} \\
Let $(\ct, \cst)$ be a local solution of (\ref{het_cat_sys}) for initial values $(\cot, \csot)$. We consider the system in the form
\begin{equation}\label{het_cat_sys_orig_form}
	\begin{array}{rclll}
			\partial_t (\ct, \cst) + \tilde{A} (\ct, \cst) &=& \tilde{F}(\ct, \cst) &\quad\matin &(0, T) \times (\Om \times \Blat),\\[0.5em]
			\tilde{B}(\ct, \cst) &=& 0 &\quad\maton &(0, T) \times \Pi,\\[0.5em]
			(\ct, \cst)\vert_{t=0} &=& (\cot, \csot) &\quad\maton &\Om \times \Blat,
		\end{array}
\end{equation}
where  $\Pi := \Bin \times \Blat \times \Bout \times \BBlat$ and where the linear part is given as
\begin{align*}
	\tilde{A} &:= \begin{pmatrix}
	U_\nabla - D_\lapl & 0\\
	-\Kad & - D_\laplBelt + \Kde
	\end{pmatrix}: D(A) \rightarrow \LpOmBlatN{p},\\[0.5em]
	D(\tilde{A}) &:= \SobSpc{2}{p}(\Om)^N \times \SobSpc{2}{p}(\Blat)^N.
\end{align*}
Note that we implicitly take the trace on $\Blat$ in the second component of $\tilde{A} (\ct, \cst)$. The nonlinearity is given as
\begin{align*}
	\tilde{F}(\ct, \cst) := \begin{pmatrix}
		0 \\ \rc(\cst)
	\end{pmatrix},
\end{align*}
and the boundary operator is given as
\begin{align*}
	\tilde{B}(\ct, \cst) := (U_\nu\ct - D_\nu\ct - \gin, -D_\nu\ct - \Kad\ct + \Kde\cst, -D_\nu\ct,-D_{\nu_\Blat}\cst) \vert_\Pi.
\end{align*}
We also set
\begin{align*}
	D_\lapl := \diag(d_1 \lapl, \dots,  d_N \lapl), \quad
	D_\laplBelt := \diag(d_1^\Blat \laplBelt , \dots, d_N^\Blat \laplBelt)
\end{align*}
and 
\begin{align*}
	U_\nabla := \diag(u \cdot \nabla , \dots, u \cdot \nabla), \quad
	U_\nu := \diag(	u \cdot \nu, \dots, u \cdot \nu)
\end{align*}
in $N$ dimensions as well as
\begin{align*}
	\Kad := \diag(k_1^\mathrm{ad}, ..., k_N^\mathrm{ad}), \quad 
	\Kde := \diag(k_1^\mathrm{de}, ...,  k_N^\mathrm{de}).
\end{align*}
Additionally,
\begin{align*}
	D_\nu := \diag(d_1 \partial_\nu, \dots, d_N \partial_\nu), \quad
	D_{\nu_\Blat} := \diag(d_1^\Blat \partial_{\nu_\Blat} , \dots,  d_N^\Blat \partial_{\nu_\Blat}).
\end{align*}
We write $\rc$ for the vector of chemical reactions $(\ric)_{i=1,\dots,N}$ and $\gin$ for the vector of inflow profiles $(\gini)_{i=1,\dots,N}$.

Now, we decompose $(\ct, \cst)$ as
\begin{align*}
	\ct = \ceq + \c, \quad \cst = \cseq + \cs
\end{align*}
s.t.\ $(\c, \cs)$ denotes the deviation from the equilibrium $(\ceq, \cseq)$. Subtracting the system for the equilibrium from the original system (\ref{het_cat_sys_orig_form}) yields
\begin{equation}\label{het_cat_sys_centered}
	\begin{array}{rclll}
		\partial_t (\c, \cs) + A (\c, \cs) &=& F(\c, \cs) &\quad\matin &(0, T) \times (\Om \times \Blat),\\[0.5em]
		(\c, \cs)\vert_{t=0} &=& (\co, \cso) &\quad\maton &\Om \times \Blat,
	\end{array}
\end{equation}
where
\begin{align*}
	A := \tilde{A}\vert_{N(B)},\quad D(A) := \{(\c, \cs) \in D(\tilde{A}) \, : \, B(\c, \cs) = 0 \}
\end{align*}
with linear boundary conditions
\begin{align*}
	B(\c, \cs) := (U_\nu\c - D_\nu\c, -D_\nu\c - \Kad\c + \Kde\cs, -D_\nu\c,-D_{\nu_\Blat}\cs) \vert_\Pi.
\end{align*}
Moreover, we set
\begin{align*}
	F(\c, \cs) := \tilde{F}(\ceq + \c, \cseq + \cs) - \tilde{F}(\ceq, \cseq)
\end{align*}
and $\co := \cot - \ceq$, $\cso := \csot - \cseq$.
System (\ref{het_cat_sys_centered}) can now be analyized w.r.t.\ its equilibrium $(0,0)$,
which is equivalent to analyzing (\ref{het_cat_sys_orig_form}) w.r.t.\ $(\ceq, \cseq)$.

Obviously, for $\rc$ given in the form as stated in the theorem we have a growth bound of type
\begin{align*}
	\abs{\rc(y)} \leqslant M (1 + \abs{y}^\gm), \qquad y \in [0, \infty)^N,
\end{align*}
for constants $M > 0$ and $\gm \in [1, \infty)$.
Therefore, we may apply \cite[Rem.\ 4.1]{maier:het_cat} to obtain the continuity of the Nemytskij operator
\begin{align*}
	\rc: \LebSpc{p\gm}(\Blat)^N \rightarrow \LebSpc{p}(\Blat)^N.
\end{align*}
Due to the fact that $\Blat$ is a manifold of dimension $m=2$ this yields the continuity of
\begin{align*}
	F: \IOp^N \times \ISp^N \rightarrow \LpOmBlatN{p},
\end{align*}
since $\ISp \hookrightarrow \LebSpc{p \gm}(\Blat)$.
Here, we use that $2 - \frac{2}{p} - \frac{2}{p} \geqslant - \frac{2}{\gm p}$ for $p \in [2, \infty)$.
Moreover, by \cite[Rem.~4.1]{maier:het_cat} we obtain for $r > 0$ that
\begin{align*}
	\norm{\LpOmBlatN{p}}{F(\c, \cs) - F(z,z^\Blat)} &\leqslant\; \norm{\LebSpc{p}(\Blat)^N}{\rc(\cseq + \cs) - \rc(z^\Blat + \cseq)} \\
	&\leqslant C(r, \cseq) \norm{\LebSpc{p\gm}(\Blat)^N}{\cs - z^\Blat}\\ 
	&\leqslant C(r, \cseq) \norm{\IOp^N  \times \ISp^N}{(\c - z, \cs - z^\Blat)}
\end{align*}
for $(\c, \cs), (z, z^\Blat) \in \BallCompl{r}{\IOp^N  \times \ISp^N}{0,0}$ and some $C(r, \cseq) > 0$, which shows that $F$ is locally Lipschitz.
Next we will consider the \FrechS derivative of the nonlinearity $F$ at $(0,0)$. First, we note that
\begin{align*}
	&\sum_{k=1}^{N}\partial_k \ric(\cseq)\cks = \sum_{k=1}^{N} (\al_i - \bt_i)\left(\kp_b \bt_k \left( \cseq \right)^{\bt - \ek} - \kp_f \al_k \left( \cseq \right)^{\al - \ek}\right)\cks.
\end{align*}
This motivates the introduction of
\begin{align*}
	\left. \begin{aligned}
	a_k &:= (\al_k - \bt_k),\;\\ 
	b_k &:= b_k (\cseq) := \left(\kp_b \bt_k \left( \cseq \right)^{\bt - \ek} - \kp_f \al_k \left( \cseq \right)^{\al - \ek}\right)
	\end{aligned}\right.	
\end{align*}
for $k = 1,\dots,N$, where $\al := (\al_1,\dots,\al_N)^T$ and $\bt := (\bt_1,\dots,\bt_N)^T$. We set $a := (a_1,\dots,a_N)^T$ and $b := b(\cseq) := (b_1,\dots,b_N)^T$ as in the theorem.
Now, we can write the derivative of the chemical reaction as
\begin{align*}
	\tilde{M} := \tilde{M}(\cseq) := a \otimes b = 
	{\begin{pmatrix}
		a_1 b_1 & \ldots & a_1 b_N \\
		\vdots & \ddots & \vdots \\
		a_N b_1 &\ldots & a_N b_N
		\end{pmatrix}}.
\end{align*}
It is not hard to see that $\dim(N(\tilde{M})) = N-1$, if $a$ and $b$ are linearly independent.
Furthermore, the spectrum $\sg(\tilde{M}) = \{ \lm_1,\dots,\lm_N\}$ is given by $\lm_1 = a^T b$ and $\lm_2 = \dots = \lm_N = 0$.
Note that $b$ and $\lm_1$ may depend on $x \in \Blat$ if the equilibrium $\cseq$ is non-constant.
Additionally, note that for a fixed $x \in \Blat$ the symmetric part $S := \frac{1}{2} (\tilde{M} + \tilde{M}^T)$ has the spectrum
\begin{align*}
\sg(S) = \left\{ \textstyle\frac{1}{2} \left( a^T b \pm \abs{a}\abs{b}\right), 0\right\},
\end{align*}
if $a$ and $b$ are linearly independent, and
\begin{align*}
\sg(S) = \left\{  \textstyle\frac{1}{2} \left( a^T b + \abs{a}\abs{b}\right), 0\right\},
\end{align*}
if $a$ and $b$ are linearly dependent, respectively.
Now, we denote by $M$ the derivative of the nonlinearity $F$ at $(0, 0)$ s.t.\ we obtain
\begin{align*}
	M(\cseq): \LpOmBlatN{p} \rightarrow \LpOmBlatN{p}, \quad \begin{pmatrix}
		\c \\ \cs
	\end{pmatrix} \mapsto \begin{pmatrix}
		0 & 0 \\ 0 & \tilde{M}
	\end{pmatrix} \begin{pmatrix}
		\c \\ \cs
	\end{pmatrix}
\end{align*}
as the $\LebSpc{p}$-realization of the multiplication operator corresponding to the matrix $\tilde{M}(\cseq)$.
Since $\cseq \in \SobSpc{2}{p}(\Blat)^N$, we obtain that $\tilde{M}$ is bounded on $\Blat$ and $M := M(\cseq) \in \mathscr{L}(\LpOmBlatN{p})$.

Finally, we have to show that $(A, F)$ satisfies appropriate estimates to be able to use the principle of linearized stability (cf. \cite[Chap.~6: Linearized Stability (S)]{pruess:maxreg} for weaker conditions than those in \cite{pruess:movint,pruess:conv}).
Since $A$ is a linear operator and does not depend on $(\c, \cs)$, it suffices to show the estimates for the nonlinear part $F$.

To this end, let $r > 0$. From \cite[Rem. 4.1]{maier:het_cat} and $M \in \mathscr{L}(\LpOmBlatN{p})$ we obtain
\begin{align*}
	&\norm{\LpOmBlatN{p}}{F(\c, \cs) - F(0,0) - M(\cseq) (\c, \cs)}\\
	&\qquad \leqslant\; \norm{\LebSpc{p}(\Blat)^N}{\rc(\cseq + \cs) - \rc(\cseq)} + \norm{\LpOmBlatN{p}}{M(\cseq) (\c, \cs)}\\ 
	&\qquad \leqslant\; C(r, \tilde{M}, \cseq) \left(\norm{\LebSpc{p\gm}(\Blat)^N}{\cs} + \norm{\LpOmBlatN{p}}{(\c, \cs)}\right)\\ 
	&\qquad \leqslant C(r, \tilde{M}, \cseq) \; \norm{\IOp^N  \times \ISp^N}{(\c, \cs)}
\end{align*}
for $(\c, \cs) \in \BallCompl{r}{\IOp^N  \times \ISp^N}{0,0}$ and a constant $C(r, \tilde{M}, \cseq) > 0$, which completes the necessary estimates for the principle of linearized stability.

\textbf{Step 2: Linearization.}\\
Using a linearization of first order of (\ref{het_cat_sys_centered}) w.r.t.\ $(0, 0)$ we obtain the system
\begin{equation}\label{het_cat_lin}
	\begin{array}{rclll}
		\partial_t (\c, \cs) + A_0 (\c, \cs) &=& G(\c, \cs), &\quad\matin &(0, T) \times (\Om \times \Blat),\\[0.5em]
		(\c, \cs)\vert_{t=0} &=& (\co, \cso) &\quad\maton &\Om \times \Blat,
	\end{array}
\end{equation}
where
\begin{align*}
	A_0 := A - M
\end{align*}
with
\begin{align*}
	D(A_0) := \left\{\left(\c, \cs\right) \in \SobSpc{2}{p}(\Om)^N \times \SobSpc{2}{p}(\Blat)^N \, : \, B(\c, \cs) = 0 \right\} = D(A).
\end{align*}
We note that if $(\c, \cs) \in D(A_0)$ is constant, then we have $(\c, \cs) = (0,0)$ due to $B(\c, \cs) = 0$ and (\Avi).
This implies that Poincar{\'e}'s inequality (cf. Remark \ref{poincare}) is at our disposal.
Moreover, we set
\begin{align*}
G(\c, \cs) := F(\c, \cs) - F(0,0) - M (\c, \cs).
\end{align*}
Based on Theorem \ref{thm:maxreg} we obtain maximal $\LebSpc{p}$-regularity for $A_0 =A - M$ by using the fact that the perturbations caused by sorption and chemical reaction are bounded in $\LpOmBlatN{p}$.

\textbf{Step 3: Characterization of the spectrum.}\\
In the following we will use the notation
\begin{align*}
	\Kab{\al}{\bt} := \pwrKad{\al}\pwrKde{\bt}, \qquad \al, \bt \geqslant 0.
\end{align*}
Note the special cases
\begin{align*}
	\Kab{0}{0} = \id, \quad \Kab{\al}{0} = \pwrKad{\al}, \quad \Kab{0}{\bt} = \pwrKde{\bt}
\end{align*}
and the fact that $\Kab{\al}{\bt}$ commutes with $D_\lapl$, $D_\laplBelt$, $U_\nabla$, $U_\nu$, $D_\nu$ and $D_{\nu_\Blat}$ for $\al, \bt \geqslant 0$.

Since $\Om$ and $\Blat$ are bounded, the operator $A_0$ has compact resolvent. So it is sufficient to analyze the eigenvalues of $A_0$ in order to characterize its spectrum.
Furthermore, due to the compact resolvent the spectrum of $A_0$ is $p$-invariant. We note that the operator $A_0$ is well-defined in the $\LebSpc{2}$-setting.
Consequently, we will determine the $\LebSpc{2}$-spectrum of $A_0$ and transfer the result to the other values $p \in [2, \infty) \setminus \{3\}$.

Let $(f_\Om, f_\Blat) \in D(A_0)$ be an eigenvector corresponding to the eigenvalue $\lm \in \sg(A_0)$. We set
\begin{align*}
	\begin{pmatrix}
	\c \\ \cs
	\end{pmatrix} := \begin{pmatrix}
		\Kab{1}{-1} & 0 \\ 0 & \Kab{1}{-1}
	\end{pmatrix} \begin{pmatrix}
		f_\Om \\ f_\Blat
	\end{pmatrix} \in D(A_0).
\end{align*}
Now, we obtain
\begin{align*}\label{resolvent_eq}
	&\Re \sclprd*{\LOmBlatN}{\lm \, \begin{pmatrix}
			\Kab{-1}{1} & 0 \\ 0 & \Kab{-1}{1}
		\end{pmatrix}\begin{pmatrix}
			\c \\ \cs
		\end{pmatrix}}
		{\begin{pmatrix}
			\id & 0 \\ 0 & \Kab{-1}{1}
		\end{pmatrix}\begin{pmatrix}
			\c \\ \cs
		\end{pmatrix}} \\[0.5em]
	= \ &\Re \sclprd*{\LOmBlatN}{A_0\, \begin{pmatrix}
			\Kab{-1}{1} & 0 \\ 0 & \Kab{-1}{1}
		\end{pmatrix}\begin{pmatrix}
			\c \\ \cs
		\end{pmatrix}}{\begin{pmatrix}
			\id & 0 \\ 0 & \Kab{-1}{1}
		\end{pmatrix}\begin{pmatrix}
			\c \\ \cs
		\end{pmatrix}}\\[0.5em]
	= \ &\Re \sclprd*{\LOmBlatN}{\begin{pmatrix}
		U_\nabla - D_\lapl & 0\\
		-\Kab{1}{0} & - D_\laplBelt + \Kab{0}{1} - \tilde{M}
		\end{pmatrix} \begin{pmatrix}
		\Kab{-1}{1} \c \\ \Kab{-1}{1} \cs
		\end{pmatrix}}{\begin{pmatrix}
		\c \\ \Kab{-1}{1} \cs
		\end{pmatrix}}\\[0.5em]
	= \ & \Re F_{\Om} + \Re F_{\Blat},
\end{align*}
where
\begin{align*}
	F_{\Om} &:= \sclprd{\LOmN}{\Kab{-1}{1} U_\nabla\c}{\c} - \sclprd{\LOmN}{\Kab{-1}{1} D_\lapl\c}{\c},\\[0.5em]
	F_{\Blat} &:= -\ \sclprd{\LBlatN}{\Kab{-1}{2} \c}{\cs} - \sclprd{\LBlatN}{\Kab{-2}{2} D_\laplBelt \cs}{\cs}\\[0.5em] & \qquad \qquad + \sclprd{\LBlatN}{\Kab{-2}{3} \cs}{\cs} - \sclprd{\LBlatN}{\tilde{M} \Kab{-1}{1} D_\laplBelt \cs}{\Kab{-1}{1}\cs}.
\end{align*}
This leads to
\begin{align*}
	\Re F_\Om= \Re \sum_{i=1}^N  \left(-\left( \left(\kia\right)^{-1} \kid \di\int_{\Om} \lapl\ci\conj{\ci} dx\right) + \left(\kia\right)^{-1} \kid \int_{\Om} (u \cdot \nabla)\ci \conj{\ci} dx\right)
\end{align*}
and we observe that
\begin{align*}
	\Re\left(\di\int_{\Om} \lapl\ci\conj{\ci} dx\right) &= \Re\left( \di\int_\BOm \dn\ci \conj{\ci} d\sg\right) - \di \int_\Om \abs{\nabla\ci}^2 dx\\[0.5em]
	&= \int_\Bin (u \cdot \nu) \abs{\ci}^2 d\sg - \int_\Blat \kia\abs{\ci}^2 d\sg + \Re\left(\int_\Blat \kid\cis\conj{\ci} d\sg \right)\\[0.5em] &\qquad \qquad - \di \int_\Om \abs{\nabla\ci}^2 dx,
\end{align*}
where we used Green's formula, the boundary conditions in (\ref{het_cat_lin}) and the form of the sorption rate as given by (\ref{sorp_func}).
Moreover, we have
\begin{align*}
	\int_{\Om} (u \cdot \nabla)\ci \conj{\ci} dx &= \frac{1}{2} \int_\BOm (u \cdot \nu) \abs{\ci}^2 d\sg\\[0.5em]
	&= \frac{1}{2} \int_\Bin (u \cdot \nu) \abs{\ci}^2 d\sg + \frac{1}{2} \int_\Bout (u \cdot \nu) \abs{\ci}^2 d\sg,
\end{align*}
where we used partial integration, the boundary conditions in (\ref{het_cat_lin}) and (\Av).
Putting these pieces together we obtain
\begin{align*}
	\Re F_{\Om} &= \normQuad{\LOmNN}{\Kab{-1/2}{1/2} D_\nabla\c} - \textstyle\frac{1}{2} \sclprd{\LebSpc{2}(\Bin)^N}{\Kab{-1}{1} U_\nu c}{c} + \textstyle\frac{1}{2} \sclprd{\LebSpc{2}(\Bout)^N}{\Kab{-1}{1} U_\nu c}{c} \\[0.5em] &\qquad \qquad + \normQuad{\LBlatN}{\Kab{0}{1/2}\c} - \Re \sclprd{\LBlatN}{\Kab{-1}{2}\c}{\cs},
\end{align*}
where
\begin{align*}
	D_\nabla := \diag\left(\sqrt{d_1} \nabla, ..., \sqrt{d_N} \nabla\right),\;\; D_{\nabla_{\Blat}} := \diag\left(\sqrt{d_1^\Blat} \nabla_\Blat, ..., \sqrt{d_N^\Blat} \nabla_\Blat\right).
\end{align*}
Finally, we observe that
\begin{align*}
	\Re F_\Blat &= \normQuad{\LBlatNN}{\Kab{-1}{1} D_{\nabla_\Blat}\cs} +  \normQuad{\LBlatN}{\Kab{-1}{3/2}\cs} \\[0.5em] & \qquad \qquad -\Re \sclprd{\LBlatN}{\Kab{-1}{2}\c}{\cs} - \Re \sclprd{\LBlatN}{\tilde{M} \Kab{-1}{1}\cs}{\Kab{-1}{1}\cs},
\end{align*}
as well as
\begin{equation}\label{spec_eq_1}
	\Re \lm \left( \normQuad{\LOmN}{\Kab{-1/2}{1/2}\c} + \normQuad{\LBlatN}{\Kab{-1}{1} \cs}\right) = \Re F_\Om + \Re F_\Blat.
\end{equation}
Since all norms appearing in (\ref{spec_eq_1}) are nonnegative and we have
\begin{align*}
	- \textstyle\frac{1}{2} \sclprd{\LebSpc{2}(\Bin)^N}{\Kab{-1}{1} U_\nu c}{c},\, \textstyle\frac{1}{2} \sclprd{\LebSpc{2}(\Bout)^N}{\Kab{-1}{1} U_\nu c}{c} \geqslant 0,
\end{align*}
due to (\Av), it remains to find appropriate estimates for the remaining terms.
Using the Cauchy-Schwarz inequality and Young's inequality we obtain
\begin{align*}
	2  \abs*{\Re \sclprd{\LBlatN}{\Kab{-1}{2}\c}{\cs}} &\leqslant 2  \abs*{\sclprd{\LBlatN}{\Kab{0}{1/2}\c}{\Kab{-1}{3/2}\cs}} \\[0.5em]
	&\leqslant \normQuad{\LBlatN}{\Kab{0}{1/2}\c} + \normQuad{\LBlatN}{\Kab{-1}{3/2}\cs}.
\end{align*}
Moreover, we have
\begin{align*}
	&\abs*{\Re \sclprd{\LBlatN}{\tilde{M}(\cseq) \Kab{-1}{1} \cs}{\Kab{-1}{1} \cs}} \\[0.5em]
	&\qquad = \abs*{\Re \sclprd{\LBlatN}{S(\cseq) \Kab{-1}{1} \cs}{\Kab{-1}{1} \cs} }\\[0.5em]
	&\qquad \leqslant\max\nolimits_{\Blat} \abs{S(\cseq)}_2 \normQuad{\LBlatN}{\Kab{-1}{1}\cs} \\[0.5em]
	&\qquad \leqslant \max\nolimits_{\Blat} \Big( \textstyle\frac{1}{2} \big| a^T b(\cseq)  \pm \abs{a}\abs{b(\cseq)} \big| \Big) \normQuad{\LBlatN}{\Kab{-1}{1} \cs}\\[0.5em]
	&\qquad \leqslant C_P \max\nolimits_{\Blat} \abs{a} \abs{b(\cseq)}\, \normQuad{\LBlatNN}{\Kab{-1}{1} D_{\nabla_\Blat} \cs},
\end{align*}
where $S$ denotes the symmetric part of $\tilde{M}$ and $C_P > 0$ denotes the \PoincS constant on $\Blat$, which does not depend on $\cs$.
Now, in order to obtain $\Re \lm \geqslant 0$, we only need to employ the condition
\begin{align*}
	\max\nolimits_{\Blat} \abs{a} \abs{b(\cseq)} \leqslant \frac{1}{C_P}.
\end{align*}

Now, assume that $\Re\lm = 0$. From (\ref{spec_eq_1}) we obtain $(\c, \cs) = 0$ such that $\lm \in \rh(A_0)$. Due to the fact that $\rh(A_0)$ is open we obtain that for every $\lm \in \C$ with $\Re\lm = 0$ there exists $\eps_\lm > 0$ such that $\Ball{\eps_\lm}{}{\lm} \subset \rh(A_0)$.

Additionally, we have maximal $\LebSpc{p}$-regularity for $A_0$. So $\mu + A_0$ is sectorial with angle of sectoriality $\ph_{\mu + A_0} < \frac{\pi}{2}$ for some $\mu \geqslant 0$ and we obtain
\begin{displaymath}
\C_{\eps} := \{ \lm \in \C \, : \, \Re\lm < \eps \} \subset \rh(A_0)
\end{displaymath}
for some $\eps > 0$. This yields $\Re \lm \geqslant \eps$ for every $\lm \in \sg(A_0)$.
An application of the principle of linearized stability (cf. \cite{pruess:conv,pruess:movint}) now yields the result. \hfill $\square$

\pagebreak
\begin{remark}
	In the case that $(\ceq, \cseq)$ is an equilibrium of chemical balance (cf. Remark \ref{eq_chem_bal}) and $a = \phi b$ for some $\phi \in \R$ the situation simplifies as follows:
	The spectrum of the symmetric part $S$ of $\tilde{M}$ consists of the eigenvalues
	\begin{align*}
		\lm_1 =  \phi |b|^2, \ \lm_2 = \dots = \lm_N = 0,
	\end{align*}
	such that we obtain stability immediately if $\phi \leqslant 0$ due to the fact that the corresponding bilinear form is negative semidefinite. Since $\cseq$ fulfills the chemical balance equations, we have
	\begin{align*}
		\cseq_i b_i &= \left(\kp_b \bt_i \left( \cseq \right)^\bt - \kp_f \al_i \left( \cseq \right)^\al \right)\\[0.5em]
		&= -(\al_i - \bt_i) \kp_f  \left(\cseq \right)^\al = -\kp_f  \left(\cseq \right)^\al a_i
	\end{align*}
	such that we indeed have $\phi \leqslant 0$ and therefore stability in $L_p$ for $p \in [2, \infty) \setminus \{ 3 \}$ holds without further conditions on $a$ and $b$.
	For $\al \neq \bt$ such an equilibrium does exist, since we can set
	\begin{align*}
		\cns{1,*} = \dots = \cns{N, *} =: \gm > 0,
	\end{align*}
	where $\gm$ is determined as
	\begin{equation*}
		\begin{array}{rcl}
			{\displaystyle \kp_b \prod_{i=1}^{N} \gm^{\bt_k} - \kp_f \prod_{i=1}^{N} \gm^{\al_k}} & = & 0 \\[1.5em]
			\Leftrightarrow \ \kp_b\gm^{|\bt|} - \kp_f \gm^{|\al|} & = & 0 \\[0.5em]
			\Leftrightarrow \ {\displaystyle \left(\frac{\kp_b}{\kp_f} \right)^{\frac{1}{|\al| - |\bt|}}} & = & \gm.
		\end{array}
	\end{equation*}
	Note that $|\alpha| - |\beta| \neq 0$, if $\alpha \neq \beta$.
\end{remark}

%%%%%%%%%%%%%%%%%%%%%%%%%%%%%%%%%%%%%%%%%%%%%%%%%%%%%%%%%%%%%%%%%%%%%%%%%%%%%
%%%%%%%%%%%%%%%%%%%%%%%%%%%%%%%%%%%%%%%%%%%%%%%%%%%%%%%%%%%%%%%%%%%%%%%%%%%%%
\section{Further results on instability}\label{sec_instab}
Next we want to find sufficient conditions, which ensure an equilibrium to be unstable.
In contrast to the situation regarding stability, we now can now drop the condition (\Avi).
\begin{theorem}
	Let $p \in [2, \infty) \setminus \{3\}$, let $T = \infty$ and let $\gini \in \Gip$ for $i = 1, \dots, N$.
	Let the sorption rates be given as 
	\begin{displaymath}
		\ris(\ci, \cis) := \kia \ci - \kid \cis, \qquad i = 1,\dots,N,
	\end{displaymath}
	and the reaction rates as
	\begin{displaymath}
		\ric(\cs) := (\al_i - \bt_i)\left(\kp_b \left(\cs\right)^\bt - \kp_f \left(\cs\right)^\al\right), \qquad i = 1,\dots,N,
	\end{displaymath}
	with $\kia, \kid > 0$, $\kp_b, \kp_f > 0$ and $\al, \bt \in (\{ 0 \} \cup [1, \infty))^N \setminus \{ 0 \}^N$.
	Assume that $(\ceq, \cseq)$ is an equi- librium of (\ref{het_cat_sys}) satisfying (\APe), (\ARe), (\AIe)
	and that there exists an eigenvector $(\c, \cs)$ of $A_0$ s.t.
	\begin{align}\label{instab:cond1}
		\sclprd{\LBlatN}{b(\cseq) \cs}{a \cs} > \abs*{\sclprd{\LOmBlatN}{A (\c, \cs)}{(\c, \cs)}},
	\end{align}
	where $A_0$, $A$, $a$ and $b(\cseq)$ are defined as in Theorem \ref{thm: stability}.
	Then the equilibrium $(\ceq, \cseq)$ is unstable in $\IOp^N \times \ISp^N$ and there exists a constant $\rh > 0$ s.t.\ for every $\eta > 0$ there exists
	\begin{align*}
		(\co, \cso) \in \BallCompl{\eta}{\IOp^N \times \ISp^N}{\ceq, \cseq},
	\end{align*}
	which in case $p > 3$ satisfies (\ref{compatibility_condition}), such that the corresponding solution $(\c, \cs)$ to (\ref{het_cat_sys}) satisfies
	\begin{align*}
		\norm{\IOp^N \times \ISp^N}{(\c(t_\eta), \cs(t_\eta)) - (\ceq, \cseq)} > \rh
	\end{align*}
	for some finite time $t_\eta > 0$.
\end{theorem}
\pagebreak
\textit{Proof:} Let $(\c, \cs)$ be an eigenvector of $A_0$ corresponding to the eigenvalue $\lm \in \sg(A_0)$ and fulfilling the assumptions.
As in the proof of Theorem \ref{thm: stability} it is sufficient to work in $\LebSpc{2}$ due to the compact resolvent of $A_0$.
A multiplication of the equations with $(\c, \cs)$ yields
\begin{align*}
	\Re \lm (\normQuad{\LOmN}{\c^{}} + \normQuad{\LBlatN}{\cs}) = \sclprd{\LOmBlatN}{(A - M)\, (\c, \cs)}{(\c, \cs)}
\end{align*}
and
\begin{align*}
	\sclprd{\LOmBlatN}{A\,  (\c, \cs)}{(\c, \cs)} \geqslant 0.
\end{align*}
Using the condition (\ref{instab:cond1}) yields $\Re \lm < 0$ such that there exists a $\lm_0 \in \sg(A_0) \cap \C_-$.
The fact that $A_0$ has compact resolvent implies that the spectrum consists of isolated eigenvalues.
Since $\mu + A_0$ is sectorial with angle of sectoriality $\ph_{\mu + A_0} < \frac{\pi}{2}$ for some $\mu \geqslant 0$, we have that $\sg(A_0) \cap \C_-$ is compact and obtain a spectral gap in $\C_-$, i.e. there exists a $\dl \in (\Re \lm_0, 0)$ such that $\sg(A_0) \cap [\dl + i\R] = \emptyset$.
Now, an application of \cite[Thm. 5.4.1]{pruess:movint} yields the result. \hfill $\square$

\begin{remark}\label{remark_instability}
	We shortly note the following facts. \\[-1.5\baselineskip]
	\begin{enumerate}
		\item Dropping the condition (\Avi) extends $D(A_0)$ to constant functions in general.
		\item Let $\lm \in \sg(A_0) \setminus \{ 0 \}$. Then there exists no constant eigenvector $(\c, \cs) \in D(A_0) \setminus \{ 0 \}$ for $\lm$. In fact, let $(\c, \cs)$ be such an eigenvector.
			Due to $B(\c, \cs) = 0$ we immediately obtain $\Kad \c = \Kde \cs$ and, therefore, $\c \neq 0$ and $\cs \neq 0$.
			In view of
		\begin{align*}
			A_0 (\c, \cs) = (A - M) (\c, \cs) = (0, \tilde{M} \cs) = \lm (\c, \cs)
		\end{align*}
		this leads to a contradiction to the assumption that $(\c, \cs)$ is an eigenvector.
		\item In general it is not clear, if an eigenvector fulfilling condition (\ref{instab:cond1}) exists.
		In particular, the conditions $b(\cseq)^T\cs, a^T \cs \neq 0$ has to be fulfilled in such a case. Observe that, due to the fact that $A_0$ is not normal in general, it is not clear if there exists a basis of $\LOmBlatN$ consisting of eigenvectors of $A_0$.
	\end{enumerate}
\end{remark}

%%%%%%%%%%%%%%%%%%%%%%%%%%%%%%%%%%%%%%%%%%%%%%%%%%%%%%%%%%%%%%%%%%%%%%%%%%%%%
%%%%%%%%%%%%%%%%%%%%%%%%%%%%%%%%%%%%%%%%%%%%%%%%%%%%%%%%%%%%%%%%%%%%%%%%%%%%%
\section{Conclusion}\label{sec_concl}
In this paper we dealt with stability and instability of a heterogeneous catalysis model in a cylindrical domain. One feature of the model is the coupling of equations in the bulk and nonlinear equations on the lateral surface of the cylinder, modeling the chemical reaction which occurs during the catalysis process.

Based on previous results regarding the maximal regularity of the linearized equations we showed a stability result in the $\Lp$-setting that indicates that the behavior of solutions near stationary points of the system is determined by the chemical reactions. In our result, stability of equilibria is given dependent on a bound on the first derivative of the chemical reaction rates. As an example we considered the equilibria of chemical balance; cf.~Remark~\ref{eq_chem_bal}.

Based on the stability analysis we extracted a sufficient condition for instability, too. It seems to be difficult to give a concrete example fulfilling these conditions for instability; cf.~also~Remark~\ref{remark_instability}.
%Several approaches to this end were made, e.g.\ considering the concrete chemical equilibria for several types of chemical reactions and a reduction to the half-space setting in order to apply the Fourier transform in two directions. In particular, the attempts to find (abstract or concrete) eigenvectors fulfilling the conditions of the instability result suffered from some hindrances, cf. Remark \ref{remark_instability} for more details.
Consequently, a detailed characterization of instability for the heterogeneous catalysis model (\ref{het_cat_sys}) is left for future considerations.

\bigskip

{\bf Acknowledgements.} \ The work of C.\ Gesse and J.\ Saal was supported
by the DFG (German Science Foundation) Grant SA 1043/3-1.
%%%%%%%%%%%%%%%%%%%%%%%%%%%%%%%%%%%%%%%%%%%%%%%%%%%%%%%%%%%%%%%%%%%%%%%

\bibliography{references}

\begin{thebibliography}{10}

\bibitem{aris:mathematical_theory}
R.~{Aris}.
\newblock {\em {The mathematical theory of diffusion and reaction in permeable
  catalysts}}, volume Volume I/II.
\newblock Claredon Press, Oxford, 1975.

\bibitem{augner:het_cat}
B.~{Augner} and D.~{Bothe}.
\newblock {Analysis of some heterogeneous catalysis models with fast sorption
  and fast surface chemistry}.
\newblock {\em Journal of Evolution Equations}, 21(3):3521--3552, 2021.

\bibitem{augner:fast_sorption}
B.~{Augner} and D.~{Bothe}.
\newblock {The fast-sorption and fast-surface-reaction limit of a heterogeneous
  catalysis model}.
\newblock {\em Discrete and Continuous Dynamical Systems - S}, 14(2):533--574,
  2021.

\bibitem{maier:het_cat}
Dieter {Bothe}, Matthias {K\"ohne}, Siegfried {Maier}, and J\"urgen {Saal}.
\newblock {Global Strong Solutions for a class of Heterogeneous Catalysis
  Models}.
\newblock {\em J. Math. Anal. Appl.}, 445(1):677--709, 2017.

\bibitem{levenspiel:chemical_reaction_eng}
O.~{Levenspiel}.
\newblock {\em {Chemical Reaction Engineering}}.
\newblock Wiley-VCH, 1998.

\bibitem{morgan:global_well_posedness}
J.~{Morgan} and B.Q. {Tang}.
\newblock Global well-posedness for volume-surface reaction-diffusion systems,
  2021.

\bibitem{pruess:movint}
J.~{Pr\"uss} and G.~{Simonett}.
\newblock {\em {Moving interfaces and quasilinear parabolic evolution
  equations}}.
\newblock Birkh\"auser, Cham, 2016.

\bibitem{pruess:conv}
J.~{Pr\"uss}, G.~{Simonett}, and R.~{Zacher}.
\newblock {On convergence of solutions to equilibria for quasilinear parabolic
  problems}.
\newblock {\em J. Differential Equations}, 246:3902--3931, 2008.

\bibitem{pruess:maxreg}
Jan {Pr\"uss}.
\newblock {Maximal Regularity for Evolution Equations in $L_p$ -Spaces}, 2002.

\bibitem{schnaubelt:dynamic_nonlinear_bc}
Roland {Schnaubelt}.
\newblock {Stable and unstable manifolds for quasilinear parabolic problems
  with fully nonlinear dynamical boundary conditions}.
\newblock {\em Advances in Differential Equations}, 22(7/8):541--592, 2017.

\bibitem{soucek:het_cat}
O.~{Sou\v{c}ek}, V.~{Orava}, J.~{M\'alek}, and D.~{Bothe}.
\newblock {A continuum model of heterogeneous catalysis: Thermodynamic
  framework for multicomponent bulk and surface phenomena coupled by sorption}.
\newblock {\em International Journal of Engineering Science}, 138:82--117,
  2019.

\bibitem{tartar:SobSpc}
Luc {Tartar}.
\newblock {\em {An Introduction to Sobolev Spaces and Interpolation Spaces}}.
\newblock Number~3 in Lecture Notes of the Unione Matematica Italiana.
  Springer, 2007.

\bibitem{white:het_cat}
M.~G. {White}.
\newblock {\em {Heterogeneous Catalysis}}.
\newblock Prentice-Hall, 1990.

\end{thebibliography}
\bibliographystyle{plain}
%%%%%%%%%%%%%%%%%%%%%%%%%%%%%%%%%%%%%%%%%%%%%%%%%%%%%%%%%%%%%%%%%%%%%%
\end{document}